\documentclass[a4paper, 12pt, titlepage]{article}
\usepackage{amsmath, amsfonts}
\usepackage{breqn}
\usepackage[affil-it]{authblk}

\def\p{\partial}
\def\be{\begin{equation}}
\def\ee{\end{equation}}
\def\x{{\bf x}}
\def\u{{\bf u}}
\def\t{\tau}
\def\s{\sigma}
\def\a{\alpha}
\def\b{\beta}

\author{A H Davison \thanks{Email address: alexander.davison@wits.ac.za; Corresponding author}}
	\affil{School of Mathematics,\\ University of the Witwatersrand,\\ P Bag 3,\\ Wits,\\ 2050,\\ South Africa }
	
\author{T Sidogi \thanks{Email address: thendos@gmail.com}}
	\affil{School of Mathematics,\\ University of the Witwatersrand,\\ P Bag 3,\\ Wits,\\ 2050,\\ South Africa}

\title{Barrier Option Pricing with Symmetries}

\begin{document}

\maketitle

\begin{abstract}
\noindent
We use Lie symmetry methods to price certain types of barrier options. Usually Lie symmetry methods cannot be used to solve the Black-Scholes
equation for options because the function defining the maturity condition for an option is not smooth. However, for barrier options, this restriction
can be accommodated and a symmetry analysis utilised to find new solutions.
\end{abstract}

\section{Introduction} \label{sec:Introduction}
The Black-Scholes equation has been used to price many financial instruments, with different boundary conditions defining different types of
instruments. These boundary conditions are not always conducive to an analytic solution of the Black-Scholes equation, and in practice numerical
methods are often used.

Symmetry methods can be used to solve difficult PDEs such as the Black-Scholes equation. However, boundary conditions can often cause problems.
Usually it is assumed that for a symmetry to admit a given boundary condition, the boundary must be invariant under the symmetry, and also the
function describing the boundary condition should be invariant. Ordinarily, since Lie symmetries are smooth, any non-smooth boundary
condition will not be admitted by a Lie symmetry, and so other methods would have to be used to solve the equation.

It has been shown by Goard \cite{Goard}  that the assumption just mentioned is overly restrictive. We will summarize this in
section~\ref{sec:The Invariant Surface Condition}. Unfortunately, the work of Goard does not overcome the problem arising from
the maturity condition of an option. As is shown in \cite{BlackScholes}, if $p$ represents time, $S$ represents the price of the 
underlying asset and $V=V(S,p)$ is the price of a European call option, the boundary condition under consideration is that
when $p=T$ the option price $V(S,p)$, with strike price $K$, should be given by
	\begin{equation}
	V(S,T) = \max\{ S - K, 0\}=\begin{cases} S-K & \mbox{ if } S\geq K\\ 0 & \mbox{ if } S< K, \label{maturity}
		\end{cases}
	\end{equation}
which has a point at which the derivative with respect to $S$ does not exist. If one considers the solution surface to
the Black-Scholes equation, the condition (\ref{maturity}) forces the surface to have a fold or a crease at the point
$(K,T)$. Symmetry techniques do not accommodate such solutions, as one of the assumptions in using symmetry techniques
is that solutions are smooth.

If, however, we introduce a second boundary condition of the form
	\[
	V(S,p) = 0 \mbox{  when  } S=g(p)
	\]
where $g(p)$ is a function such that $g(T)=K$, we avoid the problem because we need not worry about values of $S$ less than $K$ for the first
boundary condition. The understanding here is that $S=g(\t)$ represents a lower barrier; if the option price falls below this value, the option
becomes worthless. In other words, we solve the Black-Scholes equation in the region given by $0\leq p \leq T$ and $S\geq g(p)$, and outside
this region the option price is taken to be zero.

\section{Notation and Preliminaries}
Let ${\bf x}=(x^1, x^2,\ldots , x^n)\in \mathbb{R}^n$ be the independent variable, with co-ordinates $x^i$. In the case of the Black-Scholes
equation, the independent variables are $p$ and $S$, where $S$ is the price of the underlying asset, and $p$ is time, which will be transformed to
${\bf x}=(x,t)$. Let
${\bf u}=(u^1, u^2,\ldots, u^m)\in \mathbb{R}^m$ be the dependent variable, with co-ordinates $u^{\alpha}$. For the Black-Scholes equation,
there is only one dependent variable, the value of the financial instrument in question, $V$, which will be transformed to ${\bf u}=u$. We use subscripts to indicate partial differentiation:
	\[
	u^{\alpha}_x = \frac{\p u^{\alpha}}{\p x}, \quad u^{\alpha}_{xx} = \frac{\p^2 u^{\alpha}}{\p x^2},
	\quad u^{\alpha}_{xt} = \frac{\p^2 u^{\alpha}}{\p x \p t},\mbox{ etc.}
	\]
We also use the notation $\u_{(1)}$ to represent the collection of all first derivatives of $\u$, and similarly
$\u_{(2)}$, $\u_{(3)}$, etc., represent collections of higher order derivatives.

An $n$-th order partial differential equation (PDE) may then be represented as
	\begin{equation}
	F\left(\x, \u, \u_{(1)}, \u_{(2)},\ldots, \u_{(n)} \right) = 0.\label{PDE}
	\end{equation}
A symmetry $X$ of such a PDE is a vector field $X=\xi^i\frac{\p}{\p x^i} + \phi^{\alpha}\frac{\p}{\p u^{\alpha}}$ (where repeated indices indicate
summation) that leaves solutions of the PDE invariant. In practice, this means that $X^{[n]}F|_{F=0} = 0$, where $X^{[n]}$ is the $n$-th
prolongation of $X$ (this is the vector field $X$ with extra terms added to show the action of $X$ on derivatives of $\u$). There are formulae
for calculating the prolongation coefficients but the details will not be relevant in this paper. The coefficients $\xi^i$ and $\phi^{\alpha}$
are assumed to be functions of $\x$ and $\u$ only (although in principle one may consider coefficients that are also functions of derivatives
of $\u$). 

It is possible (for example, see \cite{OlverRosenau}) to use a symmetry to reduce either the order or the number of variables in a PDE, by solving the system
	\[
	\frac{d\,x^1}{\xi^1} = \cdots = \frac{d\, x^n}{\xi^n} = \frac{d\, u^1}{\phi^1} = \cdots = \frac{d\, u^m}{\phi^m}
	\]
and substituting solutions into the original equation $F=0$. For a PDE with boundary conditions, however, it is also necessary that the symmetry
satisfies the invariant surface condition.

\section{The Invariant Surface Condition} \label{sec:The Invariant Surface Condition}
We briefly outline the invariant surface condition and some of its implications in this section. For a more detailed explanation, we refer the reader
to \cite{Goard}. Here we assume that $\x=(x,t)$, $\u = (u)$, and 
$X=\xi\frac{\p}{\p x} + \tau \frac{\p}{\p t} + \phi\frac{\p}{\p u}$. A solution $u=u(x,t)$ of
(\ref{PDE}) is invariant under the symmetry $X$ if and only if the invariant surface condition holds, i.e.
	\begin{equation}
	\xi(x,t,u)\frac{\p u}{\p x} + \tau(x,t,u)\frac{\p u}{\p t} = \phi(x,t,u).\label{ISC}
	\end{equation}
	
\subsection{Initial (or Terminal) Conditions}
If we impose the condition $u(x,T) = f(x)$ then substituting this condition into (\ref{ISC}), we get
	\begin{equation}
	\xi(x,T,f(x))f'(x) + \tau(x,T,f(x))\frac{\p u(x,T)}{\p t} = \phi(x,T,f(x)).\label{ISCi}
	\end{equation}
If $X$ leaves the boundary condition invariant, then this condition is automatically satisfied. On the other hand, if $X$ does not leave the
boundary condition invariant, then we solve (\ref{PDE}) for $u_t$, substitute it into (\ref{ISCi}), and solve (\ref{ISCi}) in one of two
ways: if $X$ is given then we can solve to find the most general $f$ allowed by the symmetry; if $f$ is given, then we can find the most
general symmetry of (\ref{PDE}) that admits solutions satisfying the boundary condition.

\subsection{Boundary Conditions}
Alternatively, we may wish to impose a boundary condition of the form
	\[
	u(x,t) = G(t) \mbox{  when  } x=g(t).
	\]
In this case the invariant surface condition (\ref{ISC}) becomes
	\[
	\xi(g(t),t,G(t))u_x + \tau(g(t),t,G(t))u_t = \phi(g(t),t,G(t)),
	\]
which, after some manipulation, amounts to
	\begin{equation}
	\xi_x u_x + \xi_u u^2_x + \xi u_{xx} + \tau_x u_t + \tau_u u_x u_t + \tau u_{xt} = \phi_x + \phi_u u_x. \label{barrier}
	\end{equation}
Again, this can be solved for the boundary condition, given $X$, or solved for $X$ given the boundary condition.

\section{Barrier Options}
Any option must satisfy the Black-Scholes equation, which we write here as
	\begin{equation}
	V_{p} + \frac{1}2 \sigma^2 S^2 V_{SS} + rS V_S - r V = 0. \label{BS}
	\end{equation}
where $V$ is the price of the option, $S$ is the price of the underlying commodity, and $p$ is time.

For a barrier option there will be both a terminal condition (the option matures when $p=T$) and a boundary condition, which we call a barrier. We
assume the barrier is as stated in section~\ref{sec:Introduction}, i.e. $V(S,p) = 0$ when $S=g(p)$ (so $G(p)=0$) and $g(T)=K$.

We could find symmetries and invariant solutions of this equation directly; however there are many parameters, and to simplify calculations,
we transform (\ref{BS}) to the Heat equation. The method below is not the only method; see Gazizov and Ibragimov \cite{GazizovIbragimov}.
	
\subsection{Conversion to the Heat Equation}
We break the transformation into steps:
\begin{enumerate}
	\item Let $t = \frac{\s^2}{2} (T - p)$. This reverses time, and (\ref{BS}) becomes
		\be
		-\frac{\s^2}{2} V_t + \frac{\s^2}{2} S^2 V_{SS} + rS V_S - rV = 0. \label{step1}
		\ee
	\item Let $S = K e^x$ (we assume that $S>0$; the factor $K$ introduced here simplifies computations when an initial condition is imposed),
	then (\ref{step1}) becomes
		\be
		-\frac{\s^2}{2}V_t + \frac{\s^2}{2}V_{xx} + \left(r-\frac{\s^2}{2} \right)V_x - rV = 0. \label{step2}
		\ee
	\item Let $w = e^{\a x}V$. By choosing $\a = \frac{1}2\left(\frac{2r}{\s^2} - 1\right)$, (\ref{step2}) transforms to
		\be
		-\frac{\s^2}{2}w_t + \frac{\s^2}{2}w_{xx} -\left( \frac{\s^2\a^2}{2} + r \right)w = 0. \label{step3}
		\ee
	\item Let $y=e^{\b t} w$. Choosing $\b = \frac{1}4\left(\frac{2r}{\s^2} + 1\right)^2$ gives us
	$-\frac{\s^2}{2}y_t + \frac{\s^2}{2}y_{xx} = 0$, i.e.
		\be
		y_t = y_{xx}.
		\ee
	\item Finally, let $u = \frac{1}K y$, to make the boundary conditions simpler.
\end{enumerate}

The boundary condition $V(S,T) = \max\{S-K,0 \}$ now becomes
	\[
	u(x, 0 ) = \max\left\{e^{(\a+1)x} - e^{\a x}, 0 \right\}=\begin{cases} e^{(\a+1)x} - e^{\a x} & \mbox{ if } x\geq 0\\ 0 & \mbox{ if } x< 0.
		\end{cases}
	\]

\subsection{Finite Dimensional Symmetries}
The Lie symmetries of the Heat equation $u_t = u_{xx}$ are well-known (for example, see \cite{Olver}):
	\begin{align*}
	X_1 &= 4xt \frac{\p}{\p x}  + 4t^2\frac{\p}{\p t} + u(-2t-x^2)\frac{\p}{\p u}\\
	X_2 &= x\frac{\p}{\p x} + 2t\frac{\p}{\p t} \\
	X_3 &= 2t\frac{\p}{\p x} - ux\frac{\p}{\p u} \\
	X_4 &= \frac{\p}{\p x} \\
	X_5 &= u\frac{\p}{\p u} \\
	X_6 &= \frac{\p}{\p t} \\
	X_{\infty} &= \psi(x,t) \frac{\p}{\p u} \mbox{ where } \psi_t = \psi_{xx}.
	\end{align*}
The most general finite-dimensional symmetry of the Heat equation is a linear combination of the first six of these:
	\begin{dmath*}
	X = [4c_1xt + c_2 x + 2c_3t + c_4]\frac{\p}{\p x} + [4c_1t^2 + 2c_2t + c_6]\frac{\p}{\p t} + [c_1u(-2t-x^2) - c_3ux - c_5 u]\frac{\p}{\p u}.
	\end{dmath*}
We now wish to apply the initial condition with this symmetry to (\ref{ISCi}) to find the most general symmetry that
allows the initial condition. However, due to the point of non-smoothness of the initial condition, we only consider
values $x>0$, bearing in mind that we will also be imposing a barrier condition. The condition
(\ref{ISCi}) becomes
	\begin{dmath*}
	\left( -c_1 x^2 - c_3 x - c_5 \right)\left( e^{(\a +1)x} - e^{\a x} \right) - (c_2 x + c_4) \left( (\a + 1)e^{(\a +1)x} - \a e^{\a x} \right)
			= c_6\left( (\a + 1)^2e^{(\a+1)x} - \a^2e^{\a x} \right)
	\end{dmath*}
which we solve for $c_1,\ldots, c_6$ by comparing coefficients of $e^{(\a+1)x}$, $e^{\a x}$, etc., to get
	\begin{align*}
	c_1 & = c_2 = c_3 = 0, \\
	c_4 & = -(2\a + 1) c_6, \\
	c_5 & = (\a^2 + \a)c_6.	
	\end{align*}
Without loss of generality, we set $c_6 = 1$ so that the symmetry operator becomes
	\begin{dmath}
	X = -(2\a + 1)\frac{\p}{\p x} + \frac{\p}{\p t} - (\a^2 + \a) u \frac{\p}{\p u},
	\end{dmath}
and the system
	\[
	\frac{d\, x}{\xi} = \frac{d\, t}{\tau} = \frac{d\, u}{\phi}
	\]
becomes
	\[
	\frac{d\, x}{-2\a - 1} = \frac{d\, t}{1} = \frac{d\, u}{-\a^2 - \a}.
	\]
The left-hand pair can be solved to find
	\[
	I_1 = x + (2\a + 1)t,
	\]
and the right-hand pair can be solved to find
	\[
	u = I_2 e^{(-\a^2 - \a)t},
	\]
where $I_1$ and $I_2$ are constants of integration; the theory of differential equations tells us that these invariants are functionally dependent;
we express this as follows:
	\[
	u = h(I_1)e^{(-\a^2 - \a)t} = h(x + (2\a + 1)t)e^{(-\a^2 - \a)t}.
	\]
The heat equation $u_t = u_{xx}$ now becomes
	\[
	h''(I_1) = (2\a + 1)h'(I_1) - (\a^2 + \a)h(I_1),
	\]
which can be solved to give
	\[
	h(I_1) = A e^{(\a + 1)I_1} + B e^{\a I_1},
	\]
so that
	\[
	u = A e^{(\a + 1) x + (\a + 1)^2 t} + B e^{\a x + \a^2 t}.
	\]
Now $u(x,0) = A e^{(\a + 1)x} + B e^{\a x}$, and comparing this with the initial condition for $x\geq 0$, we see that $A=1$ and $B=-1$.

Next, we transform this back into the original variables $S$, $\t$ and $V$, to get:
	\[
	V = S - Ke^{-r(T-p)}.
	\]
$V$ does indeed solve the Black-Scholes equation, and satisfies the terminal condition for $S\geq K$; by inspection we see that the only
barrier allowed by this solution is
	\[
	S = K e^{-r(T-p)}.
	\]

\subsection{Infinite Dimensional Symmetries}
We now consider symmetries of the form
	\begin{dmath*}
	X = [4c_1xt + c_2 x + 2c_3t + c_4]\frac{\p}{\p x} + [4c_1t^2 + 2c_2t + c_6]\frac{\p}{\p t} 
			+ [c_1u(-2t-x^2) - c_3ux - c_5 u + \psi(x,t)]\frac{\p}{\p u},
	\end{dmath*}
where $\psi_t = \psi_{xx}$, with the same initial condition as before, i.e. $u(x, 0 ) = \max\left\{e^{(\a+1)x} - e^{\a x}, 0 \right\}$;
the invariant surface condition (\ref{ISCi}) can now be written as
	\begin{dmath}
	\left( -c_1 x^2 - c_3 x - c_5 \right)\left( e^{(\a +1)x} - e^{\a x} \right) + \psi - (c_2 x + c_4) \left( (\a + 1)e^{(\a +1)x} - \a e^{\a x} \right)
			= c_6\left( (\a + 1)^2e^{(\a+1)x} - \a^2e^{\a x} \right) \label{ISCii}
	\end{dmath}

We start by considering functions $\psi$ of the form $\psi = A(t)e^{\a x} + B(t) e^{(\a + 1)x}$,
and since $\psi$ has to satisfy the Heat equation, we get
	\[
	\psi = k_1 e^{\a x + \a^2 t} + k_2 e^{(\a+1)x + (\a + 1)^2 t}.
	\]
Applying this to (\ref{ISCii}) and solving for arbitrary constants in a similar way to the previous subsection, we find that
	\begin{align*}
	c_1 & = c_2 = c_3 = 0, \\
	k_1 &= -\a c_4 -c_5 - \a^2 c_6, \\
	k_2 &= (\a + 1)c_4 + c_5 + (\a + 1)^2 c_6.
	\end{align*}
This gives rise to a three-dimensional sub-algebra generated by
	\begin{align*}
	X_a & = \frac{\p}{\p x} + \left(-\a e^{\a x + \a^2 t} + (\a + 1)e^{(\a+1)x + (\a + 1)^2 t}  \right)\frac{\p}{\p u},\\
	X_b & = \left(u - e^{\a x + \a^2 t} + e^{(\a+1)x + (\a + 1)^2 t}  \right)\frac{\p}{\p u},\\
	X_c & = \frac{\p}{\p t} + \left(-\a^2 e^{\a x + \a^2 t} + (\a + 1)^2e^{(\a+1)x + (\a + 1)^2 t}  \right)\frac{\p}{\p u}.
	\end{align*}
	
\subsection{Further Solutions}
We now use a linear combination of $X_a$, $X_b$ and $X_c$. The characteristic equation that we must solve is
	\[
	\frac{dx}{c_4} = \frac{dt}{c_6} = \frac{du}{-c_5 u + \psi(x,t)},
	\]
where
	\begin{dmath*}
	\psi = \left( -\alpha  c_4-c_5-\alpha ^2 c_6 \right) e^{\alpha x + \alpha ^2 t}
				 + \left((\alpha +1)c_4 + c_5 + (\alpha + 1)^2 c_6 \right) e^{(\alpha +1) x + (\alpha +1)^2 t}.
	\end{dmath*}
We solve $\frac{dx}{c_4} = \frac{dt}{c_6}$ to get $x=\frac{c_4}{c_6} t + I_1$, i.e. $I_1 = x- \frac{c_4}{c_6} t$. The equation
$\frac{dx}{c_4}= \frac{du}{-c_5 u + \psi(x,t)}$ can also be solved in one of three possible ways, depending on the arbitrary constant
$c_4$, $c_5$ and $c_6$. We assume that $\frac{c_5}{c_4}$ does not equal $\a$ or $\a + 1$ (otherwise, after a bit of work, we arrive at a
contradiction), and find
	\begin{dmath*}
	u = I_2 e^{-\frac{c_5}{c_4}x} - \frac{c_5 + \a^2 c_6 + \a c_4}{c_4 \a + c_5} e^{\a x + \a^2 t}
						+	\frac{(\a + 1)c_4 + c_5 + (\a+1)^2 c_6}{c_4(\a + 1) + c_5} e^{(\a + 1) x + (\a + 1)^2 t}.
	\end{dmath*}
We substitute this into the heat equation, and solve for $u$ using the fact that $I_2$ must be a function of $I_1$ to find that
	\begin{dmath*}
	u = A_1 e^{\frac{\left(c_4 - \sqrt{c_4^2-4 c_5 c_6}\right) \left(c_4 t-c_6 x\right)}{2 c_4^2 c_6^2}}
		+A_2 e^{\frac{\left(c_4^2-2 c_5 c_6 -\sqrt{c_4^4-4 c_4^2 c_5 c_6}\right) \left(c_4 t-c_6 x\right)}{2 c_4 c_6^2}}
		+\frac{\alpha ^2 c_6-\alpha  c_4+c_5}{\alpha  c_4+c_5} e^{\alpha x + \alpha^2 t} 
		+ \frac{ (\alpha +1)^2 c_6+(\alpha +1) c_4-c_5 }{(\alpha +1) c_4+c_5} e^{ (\a + 1)x + (\a+1)^2 t}.
	\end{dmath*}
We require that $u(x,0) = e^{(\a+1)x} - e^{\a x}$ ($x>0$), and looking at the exponentials in our expression for $u$, we see that either $A_1=A_2=0$ 
(which leads to a contradiction), or
$c_4 = -(2 \alpha  +1)c_6$ and $c_5 = \alpha  (\alpha +1) c_6$.  We arrive at the solution
	\[
	u = e^{(\a+1)x + (\a+1)^2 t} - e^{\a x + \a^2 t},
	\]
which is the same solution that we found before.

\section{Conclusion}
We have found a solution to the Black-Scholes equation given by
	\[
	V = \begin{cases} S - Ke^{-r(T-p)} &\mbox{ if }\, 0\leq p \leq T,\, S\geq Ke^{-r(T-p)} \\ 0 &\mbox{ otherwise.}\end{cases}
	\]

At first glance, it may seem that the use of the infinite-dimensional symmetry $\psi\frac{\p}{\p u}$ does not gain any extra solutions;
however, we assumed a very specific form of $\psi$; other forms may generate other solutions. We note that although use of the
invariant surface condition was crucial to finding invariant solutions that satisfied the given boundary conditions, the invariant surface
condition is necessary, but not sufficient; the invariant surface condition restricted the choices of symmetries, but in some cases,
the boundary condition(s) restricted the choices further.

In conclusion, we have shown specifically how to use symmetry techniques to find prices of barrier options for a non-standard barrier. More generally,
we have demonstrated that PDEs with non-smooth boundary conditions can be solved using symmetry techniques if further boundary conditions are
imposed to remove any non-smooth points.


\begin{thebibliography}{10}
\bibitem[1]{Goard}
	Joanna Goard, {\em Noninvariant boundary conditions}, Applicable Analysis: An International Journal 82 (2003), pp. 473--481.
	
\bibitem[2]{BlackScholes}
	F. Black and M. Scholes, {\em The pricing of options and corporate liabilities}, The Journal of Political Economy 81 (1973), pp. 637--654.
	
\bibitem[3]{OlverRosenau}
	P. J. Olver and P. Rosenau, {\em Group Invariant Solutions of Differential Equations},
	Society for Industrial and Applied Mathematics Journal of Applied Mathematics 47 (1987), pp. 263--278.

\bibitem[4]{GazizovIbragimov}
	R. K. Gazizov and N. H. Ibragimov, {\em Lie Symmetry Analysis of Differential Equations in Finance}, Nonlinear Dynamics 17 (1998),	pp. 387--407.

\bibitem[5]{Olver}
   P. J. Olver, {\em Lie groups and differential equations}, in {\em The Concise Handbook of Algebra}, A. V. Mikhalev and G. F. Pilz (eds.),
	Kluwer Academic, Dortrecht, Netherlands (2002), pp. 92--97.
	
\end{thebibliography}
\end{document}